\documentclass[a4paper]{amsart} 

\usepackage[mathscr]{eucal}
\usepackage{amssymb} 
\usepackage{latexsym}
\usepackage{amsthm}
\usepackage{amscd}
\usepackage{amsmath}
\usepackage{amsfonts}

\theoremstyle{plain}
\newtheorem{theorem}{Theorem}[section]
\newtheorem{proposition}[theorem]{Proposition}

\newtheorem{lemma}[theorem]{Lemma}

\theoremstyle{definition}
\newtheorem{definition}[theorem]{Definition}
\newtheorem{example}[theorem]{Example}

\newcommand{\bx}{\mbox{\boldmath $x$}}
\newcommand{\by}{\mbox{\boldmath $y$}}
\def\maru#1{{\ooalign{\hfill$\scriptstyle#1$\hfill\crcr$\bigcirc$}}}
\def\daen#1{{\ooalign{\hfill$\scriptstyle#1$\hfill\crcr$\put(9,2){\oval(26,10)}$}}}
\def\ddaen#1{{\ooalign{\hfill$\scriptstyle#1$\hfill\crcr$\put(15,2){\oval(36,10)}$}}}

\title[Rectangular Schur functions]
{Rectangular Schur functions and the basic
representation of affine Lie algebras}

\author{Hiroshi Mizukawa \and Hiro-Fumi Yamada}
\address{Hiroshi Mizukawa, Department of Mathematics, Hokkaido University, Sapporo 060-0810, 
Japan}
\email{mzh@math.sci.hokudai.ac.jp\quad mzh@math.okayama-u.ac.jp}
\address{Hiro-Fumi Yamada, Department of Mathematics, Okayama University, Okayama 700-8530, 
Japan}
\email{yamada@math.okayama-u.ac.jp}

\begin{document}
\maketitle
\footnote[0]{Classification number :05E05,17B65.}
\begin{abstract}  
An expression is given for the plethysm $p_{2}\circ S_{\square}$, 
where $p_{2}$ is the power sum of degree two and $S_{\square}$
is the Schur function indexed by a rectangular partition.
The formula can be well understood from the viewpoint 
of the basic representation of the affine Lie algebra of type $A_{2}^{(2)}$
\end{abstract}

\section{Introduction}\label{intro}
The aim of this paper is to prove a formula concerning Schur functions
indexed by rectangular Young diagrams.  More precisely we will give an
 expression of the plethysm $p_2 \circ S_{\square(n,m)}$ where $p_2$ is the
power sum of degree two and $S_{\square(n,m)}$ is the Schur function indexed by
the rectangular partition $(m^n)$ (Theorem \ref{mainth}).  In fact the formula is a
 version
of one given in \cite{cr}.  Our main contribution in this paper is to
understand
their formula from the viewpoint of representations of affine Lie
algebras.
By looking at the homogeneous realization of the basic representation
of the affine Lie algebra of type $A^{(2)}_2$, we can shed another light
to that formula.  Our point of view  relies on the isomorphism
between the principal and homogeneous realizations of the basic
representation,
which was first obtained by Leidwanger \cite{leid}.  As a merit of our
understanding,
it becomes clear that the formula gives an explicit expression of a
homogeneous
polynomial $\tau$-function of a hierarchy of nonlinear differential
equations.
A similar formula for the affine Lie algebra of type $A^{(1)}_1$ is
already given in \cite{iy}.

The paper is organized as follows.  In Section \ref{s-fn} we fix some notations
concerning
the symmetric functions.  
Section \ref{coreq} is devoted to a brief
 review of division calculus of partitions, namely, cores and quotients and their 
{\it bar} analogues. In Section \ref{recta11} 
we give a simple combinatorial proof 
of a formula which is presented in \cite{iy}.
 Cores and quotients of partitions are required
in proving, while in \cite{iy} the vertex operators are used in the proof.
In Section \ref{recta22} we discuss the basic representation of the affine Lie
algebra
of type $A^{(2)}_2$.  The main result is given in Section \ref{maint}.
\section{Symmetric Functions}\label{s-fn}
We denote by $P_{n}$ the set of all partitions of $n$, $SP_{n}$ the set of all strict partitions of $n$ 
and $OP_{n}$ the set of those partitions of $n$ whose parts are odd numbers.
Let $\chi^{\lambda}_{\rho}$ be the irreducible character of 
 the
symmetric group  $S_{n}$, indexed 
by $\lambda \in P_n$ and evaluated at the conjugacy class $\rho$,   
and $\zeta^{\lambda}_{\rho}$ be the irreducible {\it negative} character of the
double cover ${\tilde S}_{n}$ (cf.
\cite{hh}), 
 indexed by $\lambda \in SP_{n}$ and evaluated at the conjugacy class 
 $\rho$. 
Here we recall symmetric functions
 of variables $\bx=(x_{1},x_{2},\cdots)$
 which are discussed in this paper.
Let $p_{r}(\bx)=\sum_{i\geq1}x_{i}^{r}$ be the power sum symmetric function for $r \geq1$.
The Schur functions are defined as follows:
\begin{align*}
S_{\lambda}(\bx)&=\sum_{\rho \in
P_{n}}z_{\rho}^{-1}\chi_{\rho}^{\lambda}p_{\rho}(\bx).
\end{align*}
For $\lambda \in SP_{n}$
define Schur's $Q$-function and $P$-function by
\begin{align*}
Q_{\lambda}(\bx)&=\sum_{\rho \in
OP_{n}}2^{(l(\lambda)+l(\rho)+\epsilon(\lambda))/2}
z_{\rho}^{-1}\zeta_{\rho}^{\lambda}p_{\rho}(\bx),\\
P_{\lambda}(\bx)&=2^{-l(\lambda)}Q_{\lambda}(\bx),
\end{align*}
where
\begin{equation*}
 \epsilon(\lambda)=
\begin{cases}
0 & \text{if $n-l(\lambda)$ is even,}\\
1&  \text{if $n-l(\lambda)$ is odd}.
	\end{cases}
\end{equation*}
Let 
$\bx=(x_{1},x_{2},\cdots)$ and
$\by=(y_{1},y_{2},\cdots)$ be variables.
We write
\begin{align*}
\bx^{r}&=(x_{1}^{r},x_{2}^{r},\cdots),\\
\bx\by&=(x_{i}y_{j};i\geq 1, j \geq 1).
\end{align*}
When $\by$ is specialized as $\by=(1,\omega,\omega^2,\cdots,\omega^{r-1},0,0,\cdots)$
 for
$\omega=\exp (2\pi\sqrt {-1}/r)$,
 we write
$$\bx\omega_{r}=\bx \by.$$
For a symmetric function $F(x)$, the plethysm $p_{r}\circ F(\bx)$ with
the $r$-th power sum $p_{r}$
is by definition \cite[p135]{mac}
$$p_{r}\circ F(\bx)=F(\bx^r).$$ 
\section{Division Calculus of Partitions}\label{coreq}
Let $r$ $(r \geq 2)$ be an integer.
For a partition $\lambda=(\lambda_{1},\cdots,\lambda_{m})$, $m$ is always supposed to be a
multiple of $r$. Put $\delta_{m}=(m-1,m-2,\ldots,1,0)$ and
$\xi=(\xi_{1},\ldots,\xi_{m})=\lambda+\delta_{m}$. We set, for $k=0,1,\cdots,r-1,$ 
$$M_{k}=\{\frac{\xi_{i}-k}{r}|\xi_{i}\equiv k {\mathrm (mod}\ r\mathrm{)}\}$$
and
$$\xi^{(k)}=(\xi_{1}^{(k)},\ldots,\xi_{m_{k}}^{(k)}),\ \xi_{i}^{(k)}\in M_{k},\ (1 \leq i \leq m_{k}),\ \xi_{1}^{(k)} > \ldots >
 \xi_{m_{k}}^{(k)}.$$
We write
$$\lambda[k]=\xi^{(k)}-\delta_{m_{k}}.$$
\begin{definition}
The collection
$$\lambda^{q(3)}=(\lambda[0],\lambda[1],\ldots,\lambda[r-1])$$
is called the {\it $r$-quotient} of the partition $\lambda$.
\end{definition}
For the strict partition $\xi=\lambda
+\delta_{m}$ we put a set 
of beads on the assigned positions as follows:
\begin{align*}
{\begin{array}{cccc}
\maru{0}&{1}&\dots&\daen{r-1}\\
r&\daen{r+1}&\dots&2r-1\\
  2r   &\daen{2r-1}&\dots& 3r-1\\
\vdots & \vdots & \vdots & \vdots\\
\end{array}}
\end{align*}
This figure represents $\xi=(\cdots,2r-1,r+1,r-1,0)$.
This beads configuration is called the $r$-abacus of $\lambda$.
A position is called a hole if it is 
not occupied by a bead.
Next we set
$C_{k}=\{rs+k|0\leq s \leq m_{k}-1\}$
and
$C=\bigcup_{k=0}^{r-1} C_{k}.$
 Let $\tilde{\xi}=(\tilde{\xi}_{1},\tilde{\xi}_{2},\ldots,\tilde{\xi}_{m})$
 be defined by $\tilde{\xi}_{i} \in C,\ \tilde{\xi}_{1}> \tilde{\xi}_{2}> \ldots >
\tilde{\xi}_{m}$.
We write $$\lambda^{c(r)}=\tilde{\xi}-\delta_{m}.$$
\begin{definition}
The partition $\lambda^{c(r)}$ is called the {\it $r$-core} of the partition
$\lambda$.
\end{definition}
\begin{example}We compute the 3-core and the 3-quotient of
$\lambda=(7,7,4,4,1)$.
We see that
\begin{eqnarray*}
&&\xi=(12,11,7,6,2,0), {\ {\rm and}}\\
&&M_{0}=\{4,2,0\},\ M_{1}=\{2\},\ M_{2}=\{3,0\}.
\end{eqnarray*}
Therefore
$$\lambda[0]=(2,1,0),\ \lambda[1]=(2),\ \lambda[2]=(2,0).$$
And we have 
\begin{eqnarray*}
&&C_{0}=\{6,3,0\},\ C_{1}=\{1\},\ C_{2}=\{5,2\},\\
&&\tilde{\xi}=(6,5,3,2,1,0).
\end{eqnarray*}
Therefore
$$\lambda^{c(2)}=(1,1).$$
\end{example}
Next we explain the {\it $r$-sign} through the example above.
Number the beads in the following two ways.
\begin{enumerate}
\item[(1)]
The natural  numbering according to the increasing order.
\item[(2)]
The 3-numbering according to the $\it layers$.
\end{enumerate}
{\center\hspace{3cm} natural numbering \hspace{1cm} 3-numbering}
\begin{align*}
{\begin{array}{ccc}
\maru{0}_{1}&1&\maru{2}_{2}\\
3&4&5\\
\maru{6}_{3}&\maru{7}_{4}&8\\
9&10&\maru{11}_{5}\\
\maru{12}_{6}&13&14\\
\vdots&\vdots&\vdots\\
\end{array}}\longrightarrow
{\begin{array}{ccc}
\maru{0}_{1}&1&\maru{2}_{3}\\
3&4&5\\
\maru{6}_{4}&\maru{7}_{2}&8\\
9&10&\maru{11}_{5}\\
\maru{12}_{6}&13&14\\
\vdots&\vdots&\vdots\\
\end{array}}
\end{align*}
When we compare the natural numbering  with  the 3-numbering we get a permutation
\begin{align*}
\sigma_{3}(\lambda)=\left( \begin{array}{cccccc}
1&2&3&4&5&6\\
1&3&4&2&5&6
\end{array} \right).
\end{align*}
The $3$-sign of $\lambda$, which is denoted by
 $\delta_{3}(\lambda)$, is defined to be the sign of
$\sigma_{3}(\lambda)$:
$$\delta_{3}(\lambda)={\rm sgn} \sigma_{3}(\lambda)=-1.$$

We now proceed to be the division calculus of strict partitions. Though 
the following arguments are valid for any positive odd integer $r$ (cf \cite{m,ol}),
we need in this paper only the case of $r=3$.
\begin{definition}
Let $\lambda=(\lambda_{1},\cdots,\lambda_{l})$
 be a strict partition. We define the double of $\lambda$ by
$$D(\lambda)=(\lambda_{1},\cdots,\lambda_{l}\ |\ \lambda_{1}-1,\cdots,\lambda_{l}-1),$$ 
in the Frobenius notation.  
\end{definition}
\begin{example}
Take $\lambda=(4,2,1)= {\begin{array}{cccc}\times&\times&\times&\times\\
\times&\times&&\\\times&&&\end{array}}$. 
The double of $\lambda$ is 
\[D(\lambda)=(5,4,4,1)={\begin{array}{ccccc}
\circ&\times&\times&\times&\times\\
\circ&\circ&\times&\times&\\
\circ&\circ&\circ&\times&\\
\circ&&&
\end{array}}.\]
\end{example}
There is a remarkable property of $D(\lambda)$ as follows.
\begin{proposition}{\rm (\cite{mac})}\label{BCQ}
Let $\lambda$ be a strict partition.
\begin{enumerate}
\item[(1)] \ There exist strict partitions $\lambda^{bc(3)}$ and $\lambda^{b}[0]$ such that
\begin{align*}
D(\lambda^{bc(3)})&=D(\lambda)^{c(3)},\\
D(\lambda^{b}[0])&=D(\lambda)[0].
\end{align*} \
\item[(2)] $D(\lambda)[2]$ is the partition conjugate to $D(\lambda)[1]$.
\end{enumerate}
\end{proposition}
\begin{definition}\label{BCQD}
\begin{enumerate}
\item[(1)]The strict partition 
$\lambda^{bc(3)}$
is called the $3$-{\it bar core} of $\lambda$.\
\item[(2)]The collection 
$$\lambda^{bq(3)}=(\lambda^{b}[0],\lambda^{b}[1])$$
is called the $3$-{\it bar quotient} of $\lambda$.
\end{enumerate}
\end{definition}
\begin{example} 
We compute the 3-bar quotient of 
$\lambda=(11,9,8,7,6,4,2)$. 
Adding 0's in the tail of $D(\lambda)$, if necessary, we always assume that the size 
of the vector $D(\lambda)$ is the multiple of $3$.
We see that
$$D(\lambda)+\delta_{12}=(23,21,20,19,18,16,14,11,9,7,2,0),$$
where $\delta_{12}=(11,10,9,8,7,6,5,4,3,2,1,0)$.
To compute  the 3-quotient of 
$D(\lambda)$, we put a set of beads on the positions 
assigned by $D(\lambda)+\delta_{12}$ as follows. 
\begin{align*}
{\begin{array}{ccc}
\maru{0}&1&\maru{2}\\
3&4&5\\
6&\maru{7}&8\\
\maru{9}&10&\maru{11}\\
12&13&\maru{14}\\
15&\maru{16}&17\\
\maru{18}&\maru{19}&\maru{20}\\
\maru{21}&22&\maru{23}\\
\vdots&\vdots&\vdots
\end{array}}
\rightarrow((4,4,2),(4,4,2),(3,3,2,2)).
\end{align*}
Read each runner from the bottom and count the 
number of holes above each bead.
Thus we have the corresponding partition
$D(\lambda)[i]$ ($i=0,1,2$).
In this case the 3-quotient of $D(\lambda)$ reads $((4,4,2),(4,4,2),(3,3,2,2))$.
This gives the 3-bar quotient of $\lambda$;
$$\lambda^{bq(3)}=((3,2),(4,4,2)).$$
Moving each bead upwards in the runner successively as far as possible, we get
\begin{align*}
{\begin{array}{ccc}
\maru{0}&\maru{1}&\maru{2}\\
\maru{3}&\maru{4}&\maru{5}\\
\maru{6}&\maru{7}&\maru{8}\\
\maru{9}&10&\maru{11}\\
12&13&\maru{14}\\
\vdots&\vdots&\vdots
\end{array}}\rightarrow(3,1)=D(2).
\end{align*}
Thus we see that $D(\lambda)^{c(3)}=(3,1)$ and  $\lambda^{bc(5)}=(2)$.

A standard exposition of 3-bar cores and quotients 
is as follows (\cite{ol}). We record the given strict partition 
$\lambda$ directly:
\begin{align*}
{\begin{array}{ccc}
0&1&\maru{2}\\
3&\maru{4}&5\\
\maru{6}&\maru{7}&\maru{8}\\
\maru{9}&10&\maru{11}\\
\vdots&\vdots&\vdots\end{array}}
\end{align*}
for $\lambda=(11,9,8,7,6,4,2)$.
This bead configuration is called the {\it 3-bar abacus} of $\lambda$.
One obtains the 3-bar quotient of $\lambda$
by the following game.
\begin{enumerate}
\item[(1)]\ The parts of the strict partition $\lambda^{b}[0]$ are 
nothing  but the levels of the beads in the 0-th runner.
\begin{align*}
{\begin{array}{c}
0\\
3\\
\large\maru{6}\\
\large\maru{9}\\
\vdots
\end{array}}\rightarrow \lambda^{b}[0]=(3,2).
\end{align*}\
\item[(2)]\ Make a sequence consisting of  0 and 1 
by the following (a)-(c).
\begin{enumerate}
 \item[(a)]\ Read the first runner from the top downwards attaching 1
 at the position with a bead, and 0 at the position without a bead.\
 \item[(b)]\ Read the second runner from the bottom (infinite away) upwards attaching 0
 at the position with a bead, and 1 at the position without a bead.\
 \item[(c)]\ Make a two-sided infinite sequence by concatenating two sequences obtained
 the above procedures (a) and (b); (a) to the right and (b) to the left.
\end{enumerate}
\begin{align*}
{\begin{array}{cc}
1_{0}&\maru{2}_{0}\\
\maru{4}_{1}&5_{1}\\
\maru{7}_{1}&\maru{8}_{0}\\
10_{0}&\maru{11}_{0}\\
13_{0}&14_{1}\\
\vdots&\vdots
\end{array}}
\rightarrow\underline{1}0010 
| 011\underline{0}\rightarrow\lambda^{b}[1]=(4,4,2).
\end{align*}
Here $\underline{1}=\cdots111$ and $\underline{0}=000\cdots$.\\
\item[(3)]\ Count the number of 0's left to each 1, 
and we have the partition $\lambda^{b}[1]$.
\end{enumerate}
The 3-bar core of $\lambda$ is given directly by
removing beads according to the following rule.
\begin{enumerate}
\item[(1)]\ Remove all beads from 0-th runner.\
\item[(2)]\ Remove the beads pairwise from the first and the second runners simultaneously. \
\item[(3)]\ Move the remaining beads upwards in the runner successively as far as possible.\
\end{enumerate}
The resulting bead configuration gives $\lambda^{bc(3)}$.
\begin{align*}
{\begin{array}{ccc}
0&1&\maru{2}\\
3&4&5\\
6&7&8\\
9&10&11\\
\vdots&\vdots&\vdots\end{array}}
\rightarrow\lambda^{bc(3)}=(2).
\end{align*}
To define the {\it {3}-bar sign} ${\bar \delta}_{3}(\lambda)$ of the strict 
partition $\lambda$, we draw the $3$-bar abacus of $\lambda$. 
Number the beads in the following two ways.
\begin{enumerate}
\item[(1)]
The natural numbering according to the increasing order.
\item[(2)]
The 3-bar numbering according to the $\it layers$.
\begin{enumerate}
\item[(a)]\ Number the beads on the 0-th runner increasing order.\
\item[(b)]\ Number the beads on the first and the second runners according to the 
increasing order of layers.
For each pair, number the bead on the second runner 
before that on the first runner.\
\item[(c)]\ Number the rest of the beads in increasing order.
\end{enumerate}
\end{enumerate}
If we take $\lambda=(11,9,8,7,6,4,2)$ as before, then the two numbering are as follows
{\center \hspace{3cm}natural numbering \hspace{0.8cm}3-bar numbering}
\begin{align*}
{\begin{array}{ccc}
0&1&\maru{2}_{1}\\
3&\maru{4}_{2}&5\\
\maru{6}_{3}&\maru{7}_{4}&\maru{8}_{5}\\
\maru{9}_{6}&10&\maru{11}_{7}\\
\vdots&\vdots&\vdots\end{array}}
\longrightarrow
{\begin{array}{ccc}
0&1&\maru{2}_{3}\\
3&\maru{4}_{4}&5\\
\maru{6}_{1}&\maru{7}_{6}&\maru{8}_{5}\\
\maru{9}_{2}&10&\maru{11}_{7}\\
\vdots&\vdots&\vdots\end{array}}.
\end{align*}
Compare the two numberings, we get a permutation
\begin{align*}
\bar{\sigma}_{3}(\lambda)&=\left(\begin{array}{ccccccc}
1&2&3&4&5&6&7\\
3&4&1&6&5&2&7
\end{array}\right).
\end{align*}
The $3$-bar sign of $\lambda$ is defined to be the sign of $\bar{\sigma}_{3}(\lambda)$:
$$\bar{\delta}_{3}(\lambda)={\rm sgn} \bar{\sigma}_{3}(\lambda)=-1.$$
\end{example}
\section{Rectangular Schur Functions and $A_{1}^{(1)}$}\label{recta11}
In this section we give a simple combinatorial derivation of a
 formula which is presented in
\cite{iy} in connection with the basic representation of affine
 Lie algebra of type
$A_{1}^{(1)}$.
Here the Schur functions are described in terms of the so called 
Sato variables:
$u_{j}=p_{j}/j(j\geq 1)$. 
We will denote them by  $S_{\lambda}(u)$. 
For a given partition we draw the
 Young diagram and fill each cell with 0 or 1 in such a
way  that $(i,j)$-cell is numbered with 0 (resp. 1) if $i+j$ is 
even (resp. odd). 

Let $\ell$ be a positive odd integer and $\Delta_{\ell}=(\ell,\ell-1,\cdots,1,0 )$ 
be the staircase partition of length $\ell$. 
The partitions $\Delta_{\ell}$ ($\ell \geq 1$, odd) are characterized as that they
are 2-cores such that ${\boxed {1}}$ can be added in each row.
Let ${\mathcal F}_1^m{(\Delta_{\ell})}$ ($0 \leq m \leq \ell$) be the set of partitions
which  are obtained by adding $m$ ${\boxed {1}}$'s to $\Delta_{\ell}$.
It is obvious that $|{\mathcal F}_1^m{(\Delta_{\ell})}|=\binom{\ell+1}{m}$. 
\begin{theorem}\label{a11}\cite{iy}
Let $\ell$ be a positive odd integer. Then 
$$\sum_{\mu \in {\mathcal F}_1^m{(\Delta_{\ell})}}=(-1 )^{|\mu[1]|}S_{\mu[0]}(u)
S_{^t \mu[1]}=S_{\square(\ell+1-m,m)}(u-v)$$
for $0 \leq m \leq \ell$.
\end{theorem}   
Before proving this theorem, we give a definition. 
A pair of partitions $(\alpha,\beta)$ is said to be complementary 
relative to the rectangle $\square(n,m)$ if $\alpha$ is contained in
 $\square(n,m)$ and the rest $\check{\beta}
=\square(n,m)-\alpha$ is the $180^{\circ}$-rotation of $\beta$:
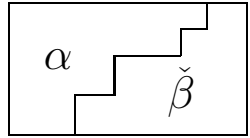
$$\begin{picture}(90,50)
\put(0,0){\line(1,0){90}}
\put(0,0){\line(0,1){50}}
\put(0,50){\line(1,0){90}}
\put(90,0){\line(0,1){50}}
\put(25,0){\line(0,1){15}}
\put(25,15){\line(1,0){15}}
\put(40,15){\line(0,1){15}}
\put(40,30){\line(1,0){25}}
\put(65,30){\line(0,1){10}}
\put(65,40){\line(1,0){10}}
\put(75,40){\line(0,1){10}}
\put(10,25){ \huge $\alpha$}
\put(60,10){\huge $\check{\beta}$}
\end{picture}$$
Note that the number of such pairs is equal to $\binom{n+m}{m}$. 
By the well-known rule for computing the Littlewood-Richardson coefficients, it is easily
seen that
\begin{align}\label{L-R}
LR_{\alpha,\beta}^{\square(n,m)}=
\begin{cases}
1&\text{if $\alpha$ and $\beta$ are complementary relative to
$\square(n,m),$}\\ 0&\text{otherwise.}
\end{cases}
\end{align}
{\it Proof of Theorem \ref{a11}}. 
It is well known \cite[p72, (5.9)]{mac} that 
$$S_{\lambda}(u-v)=\sum_{(\alpha,\beta)}(-1)^{|\beta|}LR_{\alpha,\beta}^{\lambda}
S_{\alpha}(u)S_{^t \beta}(v),$$
for any partition $\lambda$. By the above remark we have 
$$S_{\square(\ell+1-m,m)}(u-v)=\sum_{(\alpha,\beta)}(-1)^{|\beta|}S_{\alpha}(u)
S_{^t \beta}(v),$$
 where the summation runs over all pairs $(\alpha,\beta)$ of 
partitions complementary relative to the rectangle $\square(\ell+1-m,m)$.
On the other hand, by looking at the 2-abacus of $\mu \in {\mathcal
F}_1^m{(\Delta_{\ell})}$,
 we see that $\mu[0]$ and $^t \mu[1]$ are complementary relative to the rectangle
 $\square(\ell+1-m,m)$.\begin{flushright} $\square$ \end{flushright}

Theorem \ref{a11} could be understood as an expression of the weight vectors 
which are obtained by acting the group $SL_{2}$ to a maximal weight vector in
the basic representation of $A_{1}^{(1)}$. In this respect the formula gives 
an explicit description of a homogeneous polynomial $\tau$-function of the nonlinear 
Schr\"{o}dinger hierarchy, and a Virasoro singular vector as well. The reader is referred 
to \cite{iy} for the details, where the formula is proved by using the vertex 
operators. These two proofs can be transformed to each other via the isomorphism, 
given by Leidwanger \cite{leid}, between the principal and homogeneous   
realizations of the basic representation of $A_{1}^{(1)}$.
\section{Basic Representation of $A_{2}^{(2)}$}\label{recta22}
We discuss in this section the basic representation of the affine Lie algebra of type $A_{2}^{(2)}$ following
\cite{kac,kklw}. 
Here the Schur functions, Schur's $P$ and $Q$-functions are described in terms of the so called Sato variables:  
 $u_{j}=p_{j}/j\ (j \geq 1)$ for $S_{\lambda}$, $s_{j}=2p_{j}/j\ (j \geq1,\ odd)$
or $t_{j}=2p_{j}/j\
(j \geq1,\ odd)$
 for $P_{\lambda}$ and
$Q_{\lambda}$. We will denote them by $S_{\lambda}(u),\ P_{\lambda}(s),\ Q_{\lambda}(t),$ etc. 
Put $\Gamma={\mathbb C}[t_{j};\ j \geq 1,\ odd]$, whose basis is chosen as
 $\{P_{\lambda};\ \lambda \in SP_{n},\ n \in {\mathbb N} \}$.  
 Associated with the Cartan matrix \[(a_{ij})_{i,j \in \{0,1\}}=\left( \begin{array}{cc}2&-4\\ -1&2\end{array} \right),\]
 the Lie algebra ${\mathfrak g}$ of type $A_{2}^{(2)}$ is generated by $e_{i}, f_{i}, \alpha _{i}^{\vee} (i=0,1)$ and $d$ subject to the relations
 \begin{align*}
 [ \alpha^{\vee}_{i},\alpha^{\vee}_{j} ]&=0,\qquad [ \alpha^{\vee}_{i} ,e_{j}]=a_{ij}e_{j},
 \qquad [ \alpha^{\vee}_{i},f_{j}
]=-a_{ij}f_{j},\\
 [e_{i},f_{j}]&=\delta_{i,j}\alpha^{\vee}_{i},
\qquad({\rm ad} e_{i})^{1-a_{ij}}e_{j}=({\rm ad}
f_{i})^{1-a_{ij}}f_{j}=0\quad{\rm (} i\not=j{\rm )},
 \end{align*}
 and
 \begin{align*}
 [d,\alpha^{\vee}_{i}]=0,\qquad[d,e_{j}]=\delta_{j,0}e_{j},\qquad[d,f_{j}]=-\delta_{j,0}f_{j}.
 \end{align*}
 The Cartan subalgebra ${\mathfrak h}$ of ${\mathfrak g}$ is spanned by $\alpha^{\vee}_{0}, \alpha^{\vee}_{1}$ and $d$. 
 Choose the basis  $\{\alpha_{0}, \alpha_{1},\Lambda_{0}\}$ for the dual space  ${\mathfrak h}^*$ of ${\mathfrak h}$ by the pairing
 \begin{align*}
 <\alpha^{\vee}_{i},\alpha_{i}>&=a_{ij},\qquad<\alpha^{\vee}_{i},\Lambda_{0}>=\delta_{i,0},\\
 <d,\alpha_{j}>&=\delta_{0,j},\qquad<d,\Lambda_{0}>=0.
 \end{align*}
The fundamental imaginary root is $\delta=2\alpha_{0}+\alpha_{1}$.\\

The basic representation of ${\mathfrak g}$ is by definition the irreducible highest weight ${\mathfrak g}$-module with
 highest weight $\Lambda_{0}$. The weight system of the basic representation is well known:
 \[P(\Lambda_{0})=\{\Lambda_{0}-p\delta+q\alpha_{1}\ ;\ p\geq 2q^{2},\ p,q \in
\frac{1}{2}{\mathbb Z},\ p+q \in {\mathbb Z}\}.\]
 A weight $\Lambda$ on the parabola $\Lambda _{0}-2q^2\delta+q\alpha _{1}$ is said to be maximal in
the sense that $\Lambda+\delta$ 
 is no longer a weight. For any maximal weight $\Lambda$, the multiplicity of $\Lambda-n\delta\ (n \in {\mathbb N})$ is 
 known to be equal to $p(n)$, the number of partitions of $n$. A construction of the basic representation in
{\it principal} grading is 
 realized on the space $\Gamma^{(3)}={\mathbb C}[t_{j};\ j\geq 1,\ {\rm odd},\ j \not\equiv 0{\rm(mod 3)}]$ (\cite{kklw}).
 A $P$-function $P_{\lambda}(t)$ is not necessarily contained in $\Gamma^{(3)}$. However, if the strict partition $\lambda$ 
 is a 3-bar core, then $P_{\lambda}(t) \in \Gamma^{(3)}$ and in fact $P_{\lambda}(t)$ is a maximal weight vector.
 More generally we ``kill" the variables $t_{3j}\ (j \geq 1,\ odd)$ in the $P$-function $P_{\lambda}(t)$ and consider the reduced 
 $P$-function:
 \[P^{(3)}_{\lambda}(t):=P_{\lambda}(t)|_{t_{3}=t_{9}=\cdots=0} \in {\Gamma^{(3)}}.\]
 It is shown in \cite{ny} that $P_{\lambda}^{(3)}(t)$ is a weight vector for any strict partition $\lambda$, and
that
 \begin{align*}
 &\{P^{(3)}_{\lambda}(t);
\ \lambda\ {\rm is\ a\ strict\ partition\ with\ no\ part\ divisible\ by\ 3}\}\\
=&\{P^{(3)}_{\lambda}(t);\ \lambda\ {\rm is\ a\ strict\ partition\ with\ }
\lambda^{bq(3)}=(\emptyset,\lambda^{b}[1] )\}
\end{align*}
form a weight basis for $\Gamma^{(3)}$. The weight of a reduced $P$-function with a given strict partition $\lambda$ is known
 as follows. Draw the Young diagram $\lambda$ and fill each cell with 0 or 1 in such a way that, in each row the sequence (010)
  repeats from the left as long as possible. If $k_{0}$ (resp. $k_{1}$) is the number of 0's (resp. 1's) written in the Young diagram,
   then the weight of the corresponding reduced $P$-function is $\Lambda_{0}-k_{0}\alpha_{0}-k_{1}\alpha_{1}$.
A removable $i$-node ($i$=0,1) is a node ${\boxed {i}}$ of the boundary of $\lambda$ which can be removed. 
An indent $i$-node ($i$=0,1)  is a concave corner on the rim of $\lambda$ where a node   ${\boxed {i}}$ can be added. The action of 
${\mathfrak g}$ to the reduced $P$-function $P_{\lambda}^{(3)}(t)$ is described as follows:
\begin{align*}
e_{i}P_{\lambda}^{(3)}=\sum_{\mu \in {\mathcal E}_{i}^{1}(\lambda)}P_{\mu}^{(3)},
\end{align*}
where ${\mathcal E}_{i}^{1}(\lambda)$  is the set of the strict partitions 
which can be obtained by removing  a removable  $i$-node from $\lambda$, 
and
\begin{align*}
f_{i}P_{\lambda}^{(3)}=\sum_{\mu \in {\mathcal F}_{i}^{1}(\lambda)}P_{\mu}^{(3)},
\end{align*}
where ${\mathcal F}_{i}^{1}(\lambda)$ is the set of the strict partitions 
which can be obtained by adding an indent $i$-node to $\lambda$. 
For instance 
\begin{align*}
e_{0}P_{(4,3,1)}^{(3)}&=P_{(4,2,1)}^{(3)}+P_{(4,3)}^{(3)},\\
f_{1}P_{(4,3,1)}^{(3)}&=P_{(5,2,1)}^{(3)}+P_{(4,3,2)}^{(3)}.
\end{align*}
Another realization of the basic representation is known, one in the homogeneous grading.
 The isomorphism between principal and homogeneous realizations is given by Leidwanger \cite{leid}. 
 Put 
 \[{\mathcal B}={\mathbb C}[u_{j},\ s_{2j-1};\ j\geq1].\]
 Define the mapping $\Phi$ by 
 \begin{align*}
\Phi\ :\ {\Gamma} \  &\widetilde{\longrightarrow}\  {\mathcal B}\otimes{\mathbb C}[q,q^{-1}],\\
             P_{\lambda}(t) \  &\longmapsto \  2^{p(\lambda)}\bar{\delta}_{{3}}(\lambda)
P_{\lambda^{b}[0]}(s)S_{\lambda^{b}[1]}(u)\otimes
			 q^{m(\lambda)},
\end{align*}
where \[p(\lambda)=\sum_{\lambda_{i}\not\equiv 0 \pmod{3} }\left[ \frac{\lambda_{i}-1}{3}\right],\]
and $m(\lambda)$ is determined by drawing the 3-bar abacus of $\lambda$:
\begin{align*}
m(\lambda)&=({\rm number\ of\ beads\ on\ the\ first\ runner\ of\ \lambda})\\
&-({\rm number\ of\ beads\ on\ the\ second\ runner\ of\ \lambda}).
\end{align*}
For example 
\[\Phi(P_{(7,5,3,1)}(t))=8P_{(1)}(s)S_{(2,1,1)}(u)\otimes q.\]
Leidwanger \cite{leid}
 shows that $\Phi$ is indeed an isomorphism and that, if we denote by $V$ the subalgebra of ${\mathcal B}$ generated by
  $u_{2j}$ and $2^{2j-1}u_{2j-1}-s_{2j-1}\ (j \geq 1)$, then 
  \[\Phi(\Gamma^{(3)})=V\otimes{\mathbb C}[q,q^{-1}].\]
  The representation of ${\mathfrak g}$ on $V\otimes{\mathbb C}[q,q^{-1}]$, which is induced by $\Phi$, 
  is the basic representation in the homogeneous grading. In fact, if we define the degree in $V\otimes{\mathbb C}[q,q^{-1}]$ by
  \[\deg f(u,s)\otimes q^{m}=2 \deg f(u,s)+m^2,\]
then deg $\Phi(P_{\lambda}^{(3)})$ is equal to the number of 0-nodes in $\lambda$.

\section{Rectangular Schur Functions and $A_{2}^{(2)}$}\label{maint}
Let $\ell$ be a positive integer and $\Lambda_{\ell}=(3\ell-2,3\ell-5,\cdots,7,4,1)$. 
Each cell of the Young diagram of $\Lambda_{\ell}$ is supposed to be 
filled with 0 or 1 as in Section \ref{recta22}. 
Let ${\mathcal F}^{m}_{1}(\Lambda_{\ell})$ $(0 \leq m \leq \ell)$ 
be the set of the strict partitions which are obtained by adding $m$ $\boxed{1}$'s
to $\Lambda_{\ell}$. It is obvious that
 $|{\mathcal F}^{m}_{1}(\Lambda_{\ell})|=\binom{\ell}{m}$.  
We are now ready to state our main result in this paper. 
\begin{theorem}\label{mainth}
\begin{align}\label{main}
\sum_{\mu \in {\mathcal F}^{m}_{1}(\Lambda_{\ell})}
\bar{\delta}_{3}(\mu)S_{{\mu}^b[1]}=\varepsilon(\ell,m)p_{2} 
 \circ S_{\square (\ell-m,m)},
\end{align}
where
 \begin{equation*}
 \varepsilon(\ell,m)=
 \begin{cases}
(-1)^{\binom{m}{2}}\ (0 \leq m \leq \frac{\ell}{2})\\
(-1)^{\binom{\ell-m+1}{2}+(\ell-m)m}\ (\frac{\ell}{2} \leq m \leq \ell).
\end{cases}
 \end{equation*}
 \end{theorem}
It is shown in \cite{ny} that, in the principal realization of 
the basic representation of $A_{2}^{(2)}$, the $P$-functions 
$P_{\Lambda_{\ell}}(t)=P_{\Lambda_{\ell}}^{(3)}(t)$ $(\ell \geq 1)$ 
are the maximal weight vectors which allow non-zero action of $f_{1}$. 
As is explained in the previous section, we have
$$\frac{1}{m!}f_{1}^{m}P_{\Lambda_{\ell}}^{(3)}
=
\sum_{\mu \in {\mathcal
F}^{m}_{1}(\Lambda_{\ell})}P_{\mu}^{(3)}.$$ 
The left-hand side of (\ref{main}) is nothing but the image of
$\frac{1}{m!}f_{1}^{m}P_{\Lambda_{\ell}}^{(3)}$ 
under the Leidwanger isomorphism $\Phi$
to the homogeneous realization 
(dropping $q^{m(\mu)}=q^{\ell-2m}$).
Note that $p(\mu)=\binom{\ell}{2}$
for all $\mu \in {\mathcal F}_{1}^{m}(\Lambda_{\ell})$.
Therefore the formula (\ref{main}) can be thought of as 
the homogeneous realization
of the weight 
vectors which are obtained by acting the group $SL_{2}$ to a 
maximal weight vector in the basic representation of $A_{2}^{(2)}$.

The rest of this section is devoted to a proof of 
Theorem \ref{mainth}. We first recall a formula for plethysm 
which is due to Chen, Garsia and Remmel \cite{cgr}. 
\begin{proposition}\label{sxal} 
$$p_{r}\circ S_{\lambda}=\sum_{\mu}\delta_{r}(\mu) LR_{\mu[0],
\cdots,\mu[r-1]}^{\lambda}S_{\mu},$$
where the summation runs over all partitions $\mu$ with empty $r$-core. 
\end{proposition}
\begin{proof}
Recall the well-known identity \cite[p64 (4.3)]{mac} 
$$\prod_{i,j}(1-x_{i}^r y_{j}^r)^{-1}=\sum_{\mu}S_{\mu}(\bx^r)S_{\mu}(\by^r).$$
The right-hand side  can be written as 
$$\prod_{i,j}\prod_{k=0}^{r-1}(1-\omega^{k}x_{i}y_{j})^{-1}
=\sum_{\mu}S_{\mu}(\bx)S_{\mu}(\by\omega_{r}).$$
where $\omega$ is a primitive $r$-th root of 1, and $\by\omega_{r}$
denotes the $r$-inflation of the variables $\by$.
A version of Littlewood's multiple formula \cite{li2}(see also \cite{LMF})
 reads 
\begin{align*}
S_{\mu}(\by\omega_{r})=
\begin{cases}
\delta_{r}(\mu)\sum_{\nu}LR_{\mu[0],\cdots,\mu[r-1]}^{\nu}S_{\nu}(\by^{r})
&\text{if $\mu$ has empty $r$-core,}\\
0& \text{otherwise}.
\end{cases}
\end{align*}
Therefore we have 
\begin{align*}
\sum_{\mu} S_{\mu}(\bx^r)S_{\mu}{(\by^r)}=
\sum_{\mu;\mu^{c(r)}=\emptyset} \delta_{r}(\mu)S_{\mu}(\bx)
\sum_{\nu}LR_{\mu[0],\cdots,\mu[r-1]}^{\nu}S_{\nu}(\by^{r})
\end{align*}
Picking up the coefficients of $S_{\lambda}(\by^r)$, we see that 
$$S_{\lambda}(\bx^r)=
\sum_{\mu;\mu^{c(r)}=\emptyset} \delta_{r}(\mu)
LR_{\mu[0],\cdots,\mu[r-1]}^{\lambda}(\mu)S_{\mu}(\bx)$$
as desired.
\end{proof}

Our second proof also relies on Littlewood's multiple formula
 and the
orthogonality of the Schur functions. Though two proofs are 
essentially the same, they look different at a first glance. 
Readers can skip the following second proof, which is included here just as 
authors' notes.
\begin{proof}
Let $\bx^{(k)}\ (k=0,1,\cdots,r-1)$ be an infinite family of variables $\bx^{(k)}
=(x^{(k)}_{1},x^{(k)}_{2},\cdots)$. A definition of the
 Littlewood-Richardson
 coefficients is as follows:
$$S_{\lambda}(\bx^{(0)},\cdots,\bx^{(r-1)})=\sum LR_{\mu^0,\cdots
,\mu^{r-1}}^{\lambda}S_{\mu^0}(\bx^{(0)})\cdots
S_{\mu^{r-1}}(\bx^{(r-1)}),$$ 
where the summation runs over the $r$-tuples of partitions $(\mu^0,\cdots,
\mu^{r-1})$. Putting $\bx^{(0)}=\cdots=\bx^{(r-1)}=\bx^r$, 
we have 
\begin{align}\label{1}
S_{\lambda}(\bx^{r},\cdots,\bx^{r})=\sum LR_{\mu^0,\cdots
,\mu^{r-1}}^{\lambda}S_{\mu^0}(\bx^{r})\cdots
S_{\mu^{r-1}}(\bx^{r}).
\end{align} 
For an $r$-tuple $(\mu^0,\cdots,\mu^{r-1})$ of partitions, let $\mu$
 denote the single partition which has empty $r$-core and has $r$-quotient 
given by $(\mu[0],\cdots,\mu[r-1])=(\mu^0,\cdots,\mu^{r-1})$.
Then, by Littlewood's multiple formula, the right-hand side of (\ref{1})
equals
\begin{align*}
\sum_{\mu}\delta_{r}(\mu)LR_{\mu[0],\cdots,\mu[r-1]}^{\lambda}S_{\mu}(
\bx\omega_{r}),
\end{align*}
where the summation runs over all partitions $\mu$ with empty $r$-core.
The formula stated in Proposition \ref{sxal} is equivalent to 
\begin{align}\label{--2}
\sum_{\mu}LR_{\mu[0],\cdots,\mu[r-1]}^{\lambda}S_{\mu}
(\bx\omega_{r})
=
\sum_{\mu}LR_{\mu[0],\cdots,\mu[r-1]}^{\lambda}S_{\mu}
(\bx,\cdots,\bx).
\end{align}
By the Frobenius formula \cite[p114]{mac} we have
$$S_{\mu}(\bx\omega_{r})=\sum_{\rho}z_{\rho}^{-1}
\chi^{\mu}_{\rho}p_{1}(\bx\omega_{r})^{m_{1}}p_{2}(\bx\omega_{r})^{m_{2}}
\cdots,$$
where the summation runs over all partitions $\rho=(1^{m_{1}}2^{m_{2}}\cdots)$.
It is easy to see that 
\begin{align*}
p_{k}(\bx\omega_{r})=
\begin{cases}
rp_{k}(\bx) & \text{if $k \equiv 0 \pmod{r}$},\\
0 & \text{otherwise}.
\end{cases}
\end{align*}
Hence, for proving (\ref{--2}), we only need to show that 
$$\sum_{\mu}\delta_{r}(\mu)LR_{\mu[0],\cdots,\mu[r-1]}^{\lambda}
\chi^{\mu}(\rho)=0$$
for the partitions $\rho$ which are not of the form $r\sigma$.
Let $n$ be fixed. For a partition $\lambda$ of $n$ and a partition $\mu$ 
of $rn$ whose $r$-core is empty, put
$$\ell_{\lambda \mu}=\delta_{r}(\mu)LR_{\mu[0],\cdots,\mu[r-1]}^{\lambda}.$$  
 Suppose that $\rho \in P_{rn}$ is not of the form $r\sigma(\sigma \in P_{n})$.
 Then, by the orthogonality of the characters of the 
symmetric group, we see that 
\begin{align*}
0&=\sum_{\mu \in P_{rn}}\chi^{\mu}(r\sigma)\chi^{\mu}(\rho)\\
&=\sum_{\mu; \mu^{c(r)}=\emptyset}
\left(
 \sum_{\lambda \in P_{n}}\ell_{\lambda \mu}\chi^{\lambda}(\sigma)
\right)
\chi^{\mu}(\rho)\\
&=\sum_{\lambda \in P_{n}}
\left( 
\sum_{\mu; \mu^{c(r)}=\emptyset}
\ell_{\lambda \mu}  \chi^{\mu}(\rho)
\right)
\chi^{\mu}(\sigma).
\end{align*}
Here the second equality is a direct consequence of Littlewood's multiple 
formula.
Since the partition $\sigma \in P_{n}$ is arbitary, 
we can conclude that 
$$\sum_{\mu;\mu^{c(r)}=\emptyset}\ell_{\lambda \mu}\chi^{\mu}(\rho)=0$$
for any partition $\lambda \in P_{n}$ as desired.
\end{proof}

A partition $\lambda=(\lambda_{1},\cdots,\lambda_{2n})$ of $2nm$ 
is said to be $(n,m)$-{\it balanced} if $\lambda_{i}+\lambda_{2n+1-i}=2m$ 
for $1 \leq i \leq n$ (\cite[p156, Def 2.1]{cr}). 
We denote by $W(n,m)$ the set of $(n,m)$-balanced partitions.
 There is a one-to-one correspondence between the $(n,m)$-balanced partitions and pairs
 $(\alpha,\beta)$ of partitions which are complementary relative to the rectangle
$\square(n,m)$:
$$\begin{picture}(160,40)
\put(0,16.5){\line(1,0){57.5}}
\put(0,16.5){\line(0,1){25}}
\put(0,41.5){\line(1,0){82.5}}
\put(45,16.5){\line(0,1){25}}
\put(57.5,16.5){\line(0,1){7.5}}
\put(57.5,24){\line(1,0){7.5}}
\put(65,24){\line(0,1){7.5}}
\put(65,31.5){\line(1,0){12.5}}
\put(77.5,31.5){\line(0,1){5}}
\put(77.5,36.5){\line(1,0){5}}
\put(82.5,36.5){\line(0,1){5}}
\put(00,16.5){\line(0,-1){25}}
\put(00,-8.5){\line(1,0){7.5}}
\put(07.5,-8.5){\line(0,1){5}}
\put(07.5,-3.5){\line(1,0){5}}
\put(12.5,-3.5){\line(0,1){5}}
\put(12.5,1.5){\line(1,0){12.5}}
\put(25,1.5){\line(0,1){7.5}}
\put(25,9){\line(1,0){7.5}}
\put(32.5,9){\line(0,1){7.5}}
\put(55,29){$\alpha$}
\put(08,6.5){${\beta}$}

\put(120,4){\line(1,0){45}}
\put(120,4){\line(0,1){25}}
\put(120,29){\line(1,0){45}}
\put(165,4){\line(0,1){25}}
\put(132.5,4){\line(0,1){7.5}}
\put(132.5,11.5){\line(1,0){7.5}}
\put(140,11.5){\line(0,1){7.5}}
\put(140,19){\line(1,0){12.5}}
\put(152.5,19){\line(0,1){5}}
\put(152.5,24){\line(1,0){5}}
\put(157.5,24){\line(0,1){5}}
\put(125,16.5){$\alpha$}
\put(150,9){$\check{\beta}$}
\end{picture}$$
\begin{lemma}\label{W=Q}
Let $\alpha$ and $\beta$ be partitions complementary relative to the 
rectangle $\square(n,m)$. If the partition $\mu$ has empty 2-core, and 
2-quotient $(\alpha,\beta)$, then $\mu$ is $(n,m)$-balanced.
\end{lemma}
\begin{proof}
Let $\mu=(\mu_{1},\cdots,\mu_{2n})$ be a partition of $2nm$ 
whose 2-core is empty. Suppose that the 2-abacus of $\mu+\delta_{2n}$ 
can be obtained by the following.
\begin{enumerate}
\item[(1)] Put $n$ beads on the 0-th runner so that the biggest 
position at most $2(n+m-1)$.\
\item[(2)]Put $n$ beads on the first runner in exactly the same way as the 
0-th runner.\
\item[(3)] Turn over the
first
runner (together with the beads) up to the position $2(n+m)-1$.
\end{enumerate} 
Looking at the resulting 2-abacus, the 2-quotient 
of $\mu$ satisfies the condition
$$\mu[0]_{i}+\mu[1]_{n+1-i}=n\ (1 \leq i \leq n).$$
Namely $\mu[0]$ and $\mu[1]$ are complementary 
relative to the rectangle $\square(n,m)$.
Then $\mu_{i}+\mu_{2n+1-i}$ counts the holes 
up to the position $2(n+m)-1$, which is independent of $i$ and equals $2m$.
This means that $\mu$ is $(n,m)$-balanced. 
\end{proof}

Combining Lemma \ref{W=Q} with (\ref{L-R}) in 
Section \ref{recta11},
the plethysm of rectangular Schur function is rather simply expressed.
\begin{proposition}\label{carrem}\cite{cl,cr}
$$p_{2}\circ S_{\square(n,m)}=\sum_{\mu \in W(n,m)}\delta_{2}(\mu)S_{\mu}.$$
\end{proposition}
Next we need to know the relation between the two sets of partitions,
 $W(\ell-m,m)$ and ${\mathcal F}_{1}^{m}(\Lambda_{\ell})$.
\begin{lemma}\label{fpsac}
$W(\ell-m,m)=\{\mu^{b}[1]|\mu \in {\mathcal F}_{1}^{m}(\Lambda_{\ell})\}.$
\end{lemma}
\begin{proof}Since 
$$|W(\ell-m,m)|=|\{\mu^{b}[1]|\mu \in {\mathcal F}_{1}^{m}(\Lambda_{\ell})\}|,$$
we only have to show that
 $\mu^{b}[1]\in W(\ell-m,m)$ for any 
 $\mu \in {\mathcal F}^{m}_{1}(\Lambda_{\ell})$.
 Recall the relation between the 3-bar quotients and the 3-quotients:
 \[\mu^{b}[1]=D(\mu)[1].\]
First we see that the 3-abacus of $D(\Lambda_{\ell})$ is 
\begin{align*}
{\begin{array}{ccc|c}
\maru{0}&\maru{1}&2&1\\
\maru{3}&\maru{4}&5&\\
&\vdots&\\
\daen{3\ell-3}&\daen{3\ell-2}&3\ell-1&\\
3\ell&\daen{3\ell+1}&3\ell+2&\ell+1\\
3\ell+3&\daen{3\ell+4}&3\ell+5&\\
&\vdots&\\
6\ell-3&\daen{6\ell-2}&6\ell-1&2\ell
\end{array}}
\end{align*}
If we add the indent 1-node to $\Lambda_{\ell}$ at the $i$-th row, then, in the 3-abacus of
$D(\Lambda_{\ell})$, the beads  at $3i-2$ and $6\ell-3i+1$ move to $3i-1$ and
 $6\ell-3i+2$,
respectively:
\begin{align*}
{\begin{array}{ccc|c}
\maru{0}&\maru{1}&2&1\\
\maru{3}&\maru{4}&5&\\
&\vdots&\\
\daen{3i-3}&3i-2&\daen{3i-1}&i\\
&\vdots&\\
\daen{3\ell-3}&\daen{3\ell-2}&3\ell-1&\\
3\ell&\daen{3\ell+1}&3\ell+2&\ell+1\\
3\ell+3&\daen{3\ell+4}&3\ell+5&\\
&\vdots&&\\
6\ell-3i&6\ell-3i+1&\ddaen{6\ell-3i+2}&2\ell-i+1\\
&\vdots&\\
6\ell-3&\daen{6\ell-2}&6\ell-1&2\ell
\end{array}}
\end{align*}
Adding indent $1$-nodes successively, we see that, in the $3$-abacus of $D(\mu)\
 (\mu \in {\mathcal
F}^{m}_{1}(\Lambda_{\ell}))$,  the beads at $i_{1}$-th, $i_{2}$-th, $\cdots$, $i_{m}$-th
rows in the first runner  shift to the second runner as well as the beads 
at $i_{1}$-th, $i_{2}$-th, $\cdots$, $i_{m}$-th rows
from the bottom.
Then 
$D(\mu)[1]_{i}+D(\mu)[1]_{2(\ell-m)+1-i}$ counts the number of the holes of the first
runner up to the
$2\ell$-th row;
 that is $2m$. This proves that $\mu^{b}[1]=D(\mu)[1] \in W(\ell-m,m)$. 
\end{proof}

Combining Proposition \ref{carrem} and Lemma \ref{fpsac} we see so far that
\begin{align*}
\sum_{\mu \in {\mathcal F}^{m}_{1}(\Lambda_{\ell})}
{\delta}_{2}(\mu^{b}[1])S_{{\mu}^b[1]}=p_{2} 
 \circ S_{\square (\ell-m,m)}.
\end{align*}
To finish the proof of Theorem \ref{mainth} we have to look at the sign factor
 carefully. Let $\mu$ be any partition. We count the number of inversions of the
permutation $\sigma_{2}(\mu)$  corresponding to $\mu$. Form the 2-abacus of $\mu$. 
For the $i$-th bead (under natural numbering)
 let $n_{0}(i)=n_{0}$ (resp $n_{1}(i)=n_{1}$) be the number of
beads on the 0-th (resp. first) runner at the smaller position
than the bead in question. For each bead on the 0-th (resp. first)
 runner we attach a non-negative integer by the following rule:
\begin{enumerate}
\item[(1)]0 $\cdots$ if $n_{0} \geq n_{1}$
(resp. $n_{0} < n_{1}+1$),\
\item[(2)]$n_{1}-n_{0}$(resp. $n_{0}-n_{1}-1$) $\cdots$  if $n_{0}<n_{1}$
(resp. $n_{0}\geq n_{1}+1$).
\end{enumerate}
If the sum of these numbers over all beads on the 2-abacus is $q$, 
then the 2-sign $\delta_{2}(\mu)$ is equal to $(-1)^{q}$.
Here we illustrate the procedure with an example: 

\begin{align*}\mu=(9,8,8,7,6,5,5,5,5,4,3,2,2,1)\in W(7,5)\rightarrow
{\begin{array}{cc}
0&\maru{1}_{0}\\
2&\maru{3}_{0}\\
\maru{4}_{2}&5\\
\maru{6}_{1}&7\\
\maru{8}_{0}&9\\
\maru{10}_{0}&\maru{11}_{1}\\
\maru{12}_{0}&\maru{13}_{1}\\
14&\maru{15}_{0}\\
16&\maru{17}_{0}\\
18&\maru{19}_{0}\\
\maru{20}_{2}&21\\
\maru{22}_{1}&23.
\end{array}}
\end{align*}
\begin{align*}
\delta_{2}(\mu)&={\small
{\rm sgn} {\left(\begin{array}{cccccccccccccc}
1&2&3&4&5&6&7&8&9&10&11&12&13&14\\
2_{0}&4_{0}&1_{2}&3_{1}&5_{0}&7_{0}&6_{1}&9_{0}&8_{1}
&10_{0}&12_{0}&14_{0}&11_{2}&13_{1}
\end{array}\right)}}\\
&= (-1)^{2+1+1+1+2+1}=1.
\end{align*}

We now assume that $\mu$ is an element of $W(\ell-m,m)$.
Let $i^*=2(\ell-m)-i+1$ for $1\leq i \leq 2(\ell-n)$.
The bead configuration of $\mu$ has the following symmetry.
If the $i$-th bead (under natural numbering) is at the position $k$, 
then the $i^*$-th bead is at the position $2\ell-k-1$.
Thanks to this symmetry, we can count the number of inversions of
 $\sigma_{2}(\mu)$ only by looking at the first half of the beads. 
Suppose that the $i^*$-th bead ($1\leq i \leq \ell-m$) is on 
the first runner and if 
$n_{0}(i^*)>n_{1}(i^*)+1$. 
Note that this condition is equivalent to saying that the $i$-th 
bead is on the 0-th runner and $n_{0}(i)\geq n_{1}(i)$.
Then by the above rule the number $n_{0}(i^*)-n_{1}(i^*)-1$ is 
attached to that bead, which is equal to $n_{0}(i)-n_{1}(i)$.
Similar considerations for other cases lead to the following rule of 
number attachment: 
For each of the first $\ell-m$ beads we attach a non-negative integer.
Suppose that the bead is on the 0-th (resp. first) runner. 
Then attach   
\begin{enumerate}
\item[(1)$^\prime$]$n_{0}-n_{1}$(resp. $n_{1}-n_{0}+1$) $\cdots$ if $n_{0} \geq
n_{1}$(resp.
$n_{0} < n_{1}+1$),\
\item[(2)$^\prime$]$n_{1}-n_{0}$(resp. $n_{0}-n_{1}-1$) $\cdots$ if $n_{0}<n_{1}$
(resp. $n_{0}\geq n_{1}+1$).
\end{enumerate}
Note that the number attached is congruent modulo 2 to the number of holes
in the smaller positions than the bead 
in question. Hence the sum of attached numbers to the first $\ell-m$ 
beads is congruent modulo 2 to $\mu_{\ell-m+1}+\cdots+\mu_{2(\ell-m)}$. 
As a consequence we see the following lemma which is stated in \cite{cr}.  
\begin{lemma}\label{sign}
For $\mu \in W(\ell-m,m)$, 
$$\delta_{2}(\mu)=(-1)^{\mu_{\ell-m+1}+\cdots+\mu_{2(\ell-m)}}.$$
\end{lemma}
As an example, the 2-sign of the above $\mu \in W(7,5)$ is computed 
by looking at the upper half of the 2-abacus. 
\begin{align*}
{\begin{array}{cc}
0&\maru{1}_{1}\\
2&\maru{3}_{2}\\
\maru{4}_{2}&5\\
\maru{6}_{1}&7\\
\maru{8}_{0}&9\\
\maru{10}_{1}&\maru{11}_{1}\\
\vdots&\vdots\\
\end{array}}.
\end{align*}
\begin{align*}
\delta_{2}(\mu)&=(-1)^{1+2+2+1+1+1}\\
&=(-1)^{5+5+4+3+2+2+1}=1.
\end{align*}
The following lemma concludes the proof of Theorem \ref{mainth}.
\begin{lemma}
 \begin{equation*}
\delta_{2}(\mu^{b}[1])\bar{\delta}_{3}(\mu)=
 \begin{cases}
(-1)^{\binom{m}{2}}\ (0 \leq m \leq \frac{\ell}{2})\\
(-1)^{\binom{\ell-m+1}{2}+(\ell-m)m}\ (\frac{\ell}{2} \leq m \leq \ell).
\end{cases}
 \end{equation*}
\end{lemma}
\begin{proof}
Each step will be discussed 
in a parallel  way for the complementary two cases:
 (1) $m \leq \ell/2$ and (2) $m > \ell/2$.
For $\Lambda_{\ell}=(3(\ell-1)+1,\cdots ,7,4,1)$ we define the reference point
\begin{enumerate}
\item[(1)]$\stackrel{\circ}\mu=\Lambda_{\ell}+(0,\cdots,0,
\underbrace{1,0,1,0,\cdots,0,1}_{2m-1}),$
\item[(2)]$\stackrel{\circ}\mu
=\Lambda_{\ell}+(1,\cdots,1,\underbrace{0,1,0,1,0,\cdots,0,1}_{2(\ell-m)}).$
\end{enumerate}
Let $\mu$  be obtained by
$$\mu=\Lambda_{\ell}+
(0,\cdots,0,{\stackrel{i_{1}}1},0,
\cdots,0,{\stackrel{i_{2}}1},0,\cdots,0,{\stackrel{i_{m}}1},0,\cdots,0) \in 
{\mathcal F}^{m}_{1}(\Lambda_{\ell}).$$
We define a sequence
\begin{enumerate}
\item[(1)]
$\alpha(\mu)=(\ell+1-i_{m},\ell+1-i_{m-1},\cdots,\ell+1-i_{1})
-(1,3,5,\cdots,2m-1)$,
\item[(2)]\ 
$\alpha(\mu)=(\ell+1-i_{m},\ell+1-i_{m-1},\cdots,\ell+1-i_{1})-\\
(1,3,5,\cdots,2(\ell-m)-1,
2(\ell-m)+1,2(\ell-m)+2,\cdots,\ell).$
\end{enumerate}
Note that $\alpha(\mu)$ indicates the inversions encoded in $\bar{\sigma}_{3}(\mu)$. Hence
$${\bar {\delta}}_{3}(\mu)=(-1)^{\alpha(\mu)_{1}+\cdots+\alpha(\mu)_{m}}.$$
We define a sequence $\gamma(\mu)$ by
$$\gamma(\mu)_{k}=i_{m-k+1}-(m-k+1)\ (1 \leq k \leq m).$$
Since
\begin{enumerate}
\item[(1)]
$\gamma(\stackrel{\circ}\mu)=(\ell-m,\ell-m-1,\cdots,\ell-2m+1)$,
\item[(2)]
$\gamma(\stackrel{\circ}\mu)=(\ell-m,\ell-m-1,\cdots,1,0,\cdots,0)$.
\end{enumerate}
One has
$$\alpha(\mu)+\gamma(\mu)=\gamma(\stackrel{\circ}\mu).$$
Also we see that $\gamma(\mu)=(^t \mu^{b}[1]_{m+1}
,\cdots,^t \mu^{b}[1]_{\mu_{1}})$. 
From Lemma \ref{sign} we see that
$$\delta_{2}(\mu^{b}[1])=(-1)^{2m(\ell-m)-((\ell-m)m+|\gamma(\mu)|)}
=(-1)^{m(\ell-m)+|\gamma(\mu)|},$$
where $|\gamma(\mu)|=\gamma(\mu)_{1}+\cdots+\gamma(\mu)_{m}$.
Therefore one obtains 
\begin{enumerate}
\item[(1)]
$\delta_{2}(\mu^{b}[1])\bar{\delta}_{3}(\mu)=(-1)^{m(\ell-m)+|\alpha(\mu)|+|\gamma(\mu)|}
=(-1)^{\binom{m}{2}}$
\item
[(2)]
$\delta_{2}(\mu^{b}[1])\bar{\delta}_{3}(\mu)=(-1)^{\binom{\ell-m+1}{2}+(\ell-m)m}$.
\end{enumerate}
\end{proof}
\begin{example}
(1) We consider the case of $(\ell,m)=(7,3)$. 
Then we have
$$\stackrel{\circ}\mu=(19,16,14,10,8,4,2).$$
Suppose that
$\mu=(19,17,14,10,7,4,2)$.
Now we obtain the following sequences:
\begin{align*}
\alpha(\stackrel{\circ}\mu)&=(1,3,5)-(1,3,5)=(0,0,0)\\
\alpha(\mu)&=(1,5,6)-(1,3,5)=(0,2,1)\\
\gamma(\stackrel{\circ}\mu)&=(7,5,3)-(3,2,1)=(4,3,2)\\
\gamma(\mu)&=(7,3,2)-(3,2,1)=(4,1,1).
\end{align*} 
We have
$$\delta_{2}(\mu^{b}[1])\bar{\delta}_{3}(\mu)=(-1)^{\binom{3}{2}}=-1.$$
(2) We consider the case of $(\ell,m)=(7,5)$. Then we have
$$\stackrel{\circ}\mu=(20,17,14,10,8,4,2).$$
Suppose that
$\mu=(19,17,14,11,8,5,1)$.
Now we obtain the following sequences:
\begin{align*}
\alpha(\stackrel{\circ}\mu)&=(1,3,5,6,7)-(1,3,5,6,7)=(0,0,0,0,0)\\
\alpha(\mu)&=(2,3,4,5,6)-(1,3,5,6,7)=(1,0,-1,-1,-1)\\
\gamma(\stackrel{\circ}\mu)&=(7,5,3,2,1)-(5,4,3,2,1)=(2,1,0,0,0)\\
\gamma(\mu)&=(6,5,4,3,2)-(5,4,3,2,1)=(1,1,1,1,1).
\end{align*} 
We have
$$\delta_{2}(\mu^{b}[1])\bar{\delta}_{3}(\mu)=(-1)^{\binom{3}{2}+10}=-1.$$
\end{example}

\end{document}